\documentclass[12pt]{article}
\usepackage{amsthm,amstext,amscd,latexsym,amsmath,amsfonts,amssymb}
\def\0n{_{(0)}}
\def\1n{_{(1)}}
\def\2n{_{(2)}}
\def\3{_{(3)}}

\begin{document}
\title{Quantum-Classical Correspondence and Galois Extensions
\thanks{ The work is partially sponsored by Polish Committee for Scientific Research
(KBN) under Grant No 5P03B05620.}}
\author{Wladyslaw Marcinek\\
Institute of Theoretical Physics, University
of Wroc{\l}aw,\\ Pl. Maxa Borna 9, 50-204  Wroc{\l}aw,\\
Poland}
\date{}
\maketitle

\begin{abstract}
A proposal of an algebraic model for the study of relation between a quantum
environment and certain classical particle system is given. The quantum environment
is described by a category of possible quantum states, the initial particle system is
represented by an associative algebra in the category of states. The result of
particle interactions with quantum environment is described as Galois extension. A
physical applications to highly organized structures of matter are discussed.
\end{abstract}

\section{Introduction}

The study of highly organized structures of matter leads to the investigation of some
nonstandard physical particle systems and effects. The fractional quantum Hall effect
provide an example of system with well defined internal order \cite{tsg,eha,hal,zee}.
Another interesting structures appear in the so called $\frac{1}{2}$ electronic
magnetotransport anomaly \cite{jai,dst}, high temperature superconductors or laser
excitations of electrons. In these cases certain anomalous behaviour of electron have
appear. An example is also given by the so--called Lutinger liquid \cite{hal1}. The
concept of statistical--spin liquids has been studied by Byczuk and Spalek \cite{bys}.

It should be interesting to develop an algebraic approach to the unified description
of all these new structures and effects. For this goal it is natural to assume that
the whole world is divided into two parts: a classical particle system and its quantum
environment. The classical system represents an observed reality, this particles which
really exist. The quantum environment represents all quantum possibilities to become
a part of the reality in the future. The appearance of highly organized structures is
a result of certain specific interactions. We would like to construct an algebraic
model to describe these interactions. In this paper we develop a nonstandard
algebraic formalism based on the Hopf algebra theory and Galois extensions. Our
construction is described in two steps. The first step is that the initial article is
transform under interaction into composite systems consisting quasiparticles and
quanta. Such systems represent possible results of interactions
\cite{qstat,com,gqsd,gsd}. As the second step we describe an algebra of realizations
of quantum possibilities. This step is connected with construction of an algebra
extension and with the 'decision' which possibility can be realized and which one --
not. The problem is how such 'decision' is done. We use the concept of quantum
commutativity and generalized Pauli exclusion principle \cite{qsym} for the solution
of this problem. The relation between classical and quantum theory considered
previously by some authors, see \cite{hor,haa,blj} for instance, is in our opinion
not satisfactory. Our approach is based on previously developed concept of particle
systems with generalized statistics and quantum symmetries
\cite{qstat,com,gqsd,gsd,gsi}.

\section{Fundamental assumptions on the model}

Let us consider a system of charged particles interacting with an external quantum
environment. We assume that every charge is equipped with ability to absorption and
emission of quanta of certain nature. A system which contains a charge and certain
number of quanta as a result of interaction with the quantum environment is said to be
a {\it dressed particle}. A particle dressed with a single quantum is a fictitious
particle called {\it quasiparticle}. Our model is based on the assumption that every
charged particle is transform under interaction into a composite system consisting
quasiparticles and quanta. This system represent possible result of interactions.
Note that the process of absorption of quanta by a charged particle can be desribed
as creation of quasiparticles and emission annihilate quasiparticles.

Our constructions are based on the following assumptions. We assume that the quantum
environment is represented by a tensor category $\mathcal{C}$ representing possible
quantum states of the particle system. An unital and associative algebra
$\mathcal{A}$ in the category $\mathcal{C}$ represents the classical states of the
system as a part of reality. Quanta are characterized by a given finitely generated
coquasitriangular Hopf algebra $H$. Quasiparticles are described by a new algebra
which is an extension $\mathcal{A}^{\it ext}$ of $\mathcal{A}$. Interactions are
described by right action and right coaction of $\mathcal{H}$ on the algebra
$\mathcal{A}^{\it ext}$. Obviously the right action and right coaction of $H$ on
$\mathcal{A}^{\it ext}$ must be compatible. For this goal we assume that
$\mathcal{A}^{\it ext}$ is a right $\mathcal{H}$--Hopf module. It is also natural to
assume that the algebra $\mathcal{A}$ is invariant and coinvariant with respect to
the action and coaction of $H$, respectively, i. e.
\begin{equation}
\mathcal{A}\equiv(\mathcal{A}^{{\it ext}})^{H}\equiv(\mathcal{A}^{\it ext})^{\it coH},
\end{equation}
where $(\mathcal{A}^{\it ext})^{H}$ is the set of $H$--invariants, and
$(\mathcal{A}^{\it ext})^{\it coH}$ -- $H$--coinvariants. We can consider a left
action and left coaction similarly. Then we obtain Hopf bimodule. We assume that the
right (co-)action can be transform into the left by quntum commutativity.

Composite system of two quasiparticles are described as a tensoor product
$\mathcal{A}^{\it ext}\;_{\mathcal{A}}\otimes\mathcal{A}^{\it ext}$ over
$\mathcal{A}$. The process of creation of quasiparticle as a tensor product
\begin{equation}
\otimes:\mathcal{C}\times\mathcal{C}\longrightarrow\mathcal{C}
\end{equation}
sending $\mathcal{A}^{\it ext}\times\mathcal{A}^{\it ext}$ into $\mathcal{A}^{\it
ext}\;_{\mathcal{A}}\otimes\mathcal{A}^{\it ext}$. On the other hand a composite
system of of quasiparticles and quanta are described as a tensor product
$\mathcal{A}^{\it ext}\otimes\mathcal{H}$. This product represents possible quantum
configurations as a result of quantum absorption process. Our assumption leads to the
mapping
\begin{equation}
\beta :\mathcal{A}^{\it ext}\;_{\mathcal{A}}\otimes\mathcal{A}^{\it
ext}\rightarrow\mathcal{A}^{\it ext}\otimes \mathcal{H}
\end{equation}
and the assumption that $\mathcal{A}^{\it ext}$ is the $\mathcal{H}$-Galois extension
of $\mathcal{A}$. In this way $\mathcal{A}^{\it ext}$ should be at the same time Hopf
bimodule and Galois extension. If the extended algebra $\mathcal{A}^{\it ext}$ is also
a Hopf bimodule with a nontrivial quantum commutative multiplication then we say that
$\mathcal{A}^{\it ext}$ is well defined algebra of realizations. The opposite case
when the mutiplication is degenerated is connected with the generalized Pauli
exclusion principle.

\section{Hopf modules and bimodules}

Let us assume that there is a finite Hopf algebra $H = H(m, \eta, \triangle ,
\epsilon, S)$, equipped with the multiplication $m$, the unit $\eta$, the
comultiplication $\triangle$, the counit $\epsilon$ and the antipode $S$. We use the
following notation for the coproduct in $H$: if $h \in H$, then $\triangle (h) :=
\Sigma h_{\1n} \otimes h_{\2n} \in H \otimes H$. We assume that $H$ is
coquasitriangular Hopf algebra (CQTHA). This means that $H$ equipped with a bilinear
form $b : H \otimes H \rightarrow {\bf k}$ such that
\begin{equation}
\begin{array}{l}
\Sigma b(h_{\1n}, k_{\2n}) k_{\2n} h_{\2n}
= \Sigma h_{\1n} k_{\1n} b(h_{\2n}, k_{\2n}),\\
b(h, kl) = \Sigma b(h_{\1n}, k) b(h_{\2n}, l),\\
b(hk, l) = \Sigma b(h, l_{\2n}) b(k, l_{\1n})
\end{array}
\end{equation}
for every $h, k, l \in H$. If such bilinear form $b$ exists for a given Hopf algebra
$H$, then we say that there is a {\it coquasitriangular structure} on $H$.

In this paper we assume that the Hopf algebra $H$ is a group algebra ${\bf k}G$, where
$G$ is an Abelian group and ${\bf k}$ is the field of complex numbers. The group
algebra $H := {\bf k} G$ is a Hopf algebra for which the comultiplication, the
counit, and the antipode are given by the formulae
$$
\begin{array}{cccc} \triangle (g) := g \otimes g,&\eta(g)
:= 1,&S(g) := g^{-1}&\mbox{for} \ g \in G.
\end{array}
$$
respectively. If $G$ is an Abelian group, then the coquasitriangular structure on $H
= {\bf k} G$ is given by a commutation factor $b : G\otimes G\rightarrow{\bf
k}\setminus\{0\}$ on $G$, \cite{sch,mon,qsym}.

Let us briefly recall the notion of Hopf modules, see \cite{mon} for details. If $M$
is a right $H$--comodule with coaction $\delta : M \rightarrow M \otimes H$ with
respect to a given Hopf algebra $H$, then the set $M^{coH}$ of $H$-coinvarints is
defined by the formula
\begin{equation}
M^{coH} := \{m \in M: \delta(m) = m \otimes 1\}.
\end{equation}
This means that $M^{coH}$ is trivial right $H$-comodule. If $E$ is a trivial right
$H$-comodule, then $E \otimes H$ is a nontrivial right $H$-comodule. The comodule map
$\delta : E\otimes H \rightarrow E \otimes H \otimes H$ is given by the relation
\begin{equation}
\begin{array}{c}
\delta(x \otimes h) := \Sigma \ x \otimes \Delta(h), \label{dco}
\end{array}
\end{equation}
where $x \in E, h \in H$. The inclusion of the notion of right $H$-modules leads to
the concept of Hopf modules \cite{mon}.

A right $H$-Hopf module is a $k$-linear space $M$ such that\\
(i) there is a right $H$-module action $\lhd : M \otimes H \rightarrow M$,\\
(ii) there is a right $H$-comodule map $\delta : M \rightarrow M \otimes H$, \\
(iii) $\delta$ is a right $H$-module map, this means that we have the relation
\begin{equation}
\Sigma \ (m \lhd h)\0n \otimes (m \lhd h)\1n = \Sigma \ m\0n \lhd h\1n \otimes m\1n \
h\2n,
\end{equation}
where $m \in M, h \in H, \delta (m) = \Sigma \ m\0n \otimes m\1n$, and $\triangle (h)
= \Sigma \ h\1n \otimes h\2n$.

Let $E$ be a finite--dimensional vector space equipped with a vector space basis
$\{x^a : a=1,\ldots,n\}$ and $H$ be a finite Hopf algebra equipped with generators
$\{h^i : i=1,\ldots,N\}$. It is obvious that $E$ is a Hopf module equipped with a
trivial right $H$-module and comodule structures
\begin{equation}
x^a\lhd h^i=\epsilon(h^i)x^a,\quad\delta(x^a)=x^a\otimes 1.
\end{equation}
Let us consider the tensor product $M=E\otimes H$ and define the right Hopf module
coaction
\begin{equation}
\begin{array}{c}
\delta(x^a\otimes h^i) := x^a\otimes \Delta(h^i) .\label{cod}
\end{array}
\end{equation}
We use the following notation
\begin{equation}
x^a_{h^i} := x^a\otimes h^i \quad\mbox{(no sum)}.
\end{equation}
In this way the right coaction (\ref{dco}) is given by the relation
\begin{equation}
\begin{array}{c}
\delta(x^a_{h^i}):=\Sigma \ x^a_{h^i\1n} \otimes h^i\2n , \label{hco}
\end{array}
\end{equation}
where $\triangle (h^i) = \Sigma \  h^i\1n \otimes h^i\2n$. For the right action $\lhd
: M \otimes H \rightarrow M$ we obtain
\begin{equation}
x^a_{h^i} \lhd h^j = x^a_{h^i h^j}.
\end{equation}
It this way the right $H$-Hopf module $M$ is a product $M=E\otimes H$, where $E\equiv
M^{coH}$ is the trivial right $H$-module.

The concept of left Hopf $H$--modules can be introduced in a similar way. If $M$ is a
right and left $H$--Hopf module, then it is said to be a $H$--Hopf bimodule. If $H$
is CQTHA, then every right $H$--Hopf bimodule is $H$--Hopf bimodule. In this case the
right and left coactions coincide (up to a constant). The left action $\rhd :
H\otimes M \rightarrow M$ is given by
\begin{equation}
h^j \rhd x_{h^i} = b(h^i , h^j) x_{h^i h^j}.
\end{equation}

If $H = {\bf k} G$, then every right $H$-Hopf module $M$ is a $G$-graded space $M =
\oplus_{g\in G}M_g$, such that $M_g := E \otimes g$, where $g \in G$ and $E$ is an
arbitrary linear space with trivial action of $G$. In this case we have
\begin{equation}
x^a_g\lhd h = x^a_{gh}, \;\;\; \delta (x^a_g) = x^a_g \otimes g,
\end{equation}
where $x^a_g \in M_g$, $x^a_g := x^a\otimes g$, $x^a\in E$, $g, h\in G$.

Let $\mathcal{A}$ be an unital and associative algebra and $H$ be a finite Hopf
algebra. If $\mathcal{A}$ is a right $H$--comodule such that the multiplication map
$m : \mathcal{A}\otimes\mathcal{A}\rightarrow\mathcal{A}$ and the unit one $\eta :
{\bf k}\rightarrow\mathcal{A}$ are $H$--comodule maps, then we say that it is a right
$H$--comodule algebra.

The algebra $\mathcal{A}$ is said to be {\it quantum commutative} with respect to the
coaction of $(H, b)$ if an only if we have the relation
\begin{equation}
\begin{array}{c}
a \ b = \Sigma \ b(a_{(1)}, b_{(1)}) \ b_{(0)} \ a_{(0)} ,
\end{array}
\end{equation}
where $\rho (a) = \Sigma a_{(0)}\otimes a_{(1)}\in\mathcal{A} \otimes H$, and $\rho
(b)=\Sigma b_{(0)}\otimes b_{(1)} \in \mathcal{A} \otimes H$ for every $a, b \in
\mathcal{A}$, see \cite{cowe}. The Hopf algebra $H$ is said to be a {\it quantum
symmetry} for $\mathcal{A}$.

An algebra extension $\mathcal{A}^{\it ext}$ of $\mathcal{A}$ such that it is a right
$\mathcal{H}$--comodule algebra and $\mathcal{A}$ is its coinvariant subalgebra
\begin{equation}
\begin{array}{c}
{\mathcal A}\equiv({\mathcal A}^{\it ext})^{\it co{\mathcal H}} :=\{a \in{\mathcal
A}^{\it ext}:\delta(a)=a\otimes 1\}
\end{array}
\end{equation}
is said to be $\mathcal{H}$--extension. If in addition the map $\beta
:\mathcal{A}^{\it ext}\;_{\mathcal{A}}\otimes\mathcal{A}^{\it
ext}\rightarrow\mathcal{A}^{\it ext}\otimes \mathcal{H}$ defined by
\begin{equation}
\begin{array}{c}
\beta (a\;_{\mathcal{A}}\otimes b) := (a \otimes 1)\delta (b). \label{gal}
\end{array}
\end{equation}
is bijective then the $\mathcal{H}$--extension is said to be Galois.

If $\mathcal{A}^{\it ext}$ is $\mathcal{H}$-Galois extension, then there is a
bijection
\begin{equation}
\begin{array}{c}
\beta^n:\underbrace{{\mathcal A}^{\it ext}\;_{{\mathcal A}}\otimes\cdots\;_{{\mathcal
A}}\otimes{\mathcal A}^{\it ext}}_{n+1}\rightarrow{\mathcal A}^{\it ext}
\otimes\underbrace{{\mathcal H}\otimes\cdots\otimes{\mathcal H}}_{n}
\end{array}
\end{equation}
is given by
\begin{equation}
\begin{array}{c}
\beta^n := (\beta\otimes id)\circ\cdots\circ(id\otimes\beta\otimes
id)\circ(id\otimes\beta).
\end{array}
\end{equation}
It is interesting that for the group algebra $H \equiv {\bf k} G$ the extension
$\mathcal{A}^{\it ext}$ of $\mathcal{A}$ is Galois if and only if it is strongly
$G$--graded algebra, i e.
\begin{equation}
\begin{array}{ccc}
{\mathcal{A}}^{\it ext}=\oplus_{g\in G}{\mathcal{A}}^{\it ext}_g,&{\mathcal{A}}^{\it
ext}_g{\mathcal{A}}^{\it ext}_h = {\mathcal{A}}^{\it ext}_{gh},& {\mathcal{A}}^{\it
ext}_e\equiv{\mathcal{A}},
\end{array}
\end{equation}
$e$ is the neutral element of $G$. If we assume that there is a trivial action and
coaction of $H \equiv {\bf k} G$ on ${\mathcal{A}}^{\it ext}$, then
\begin{equation}
{\mathcal{A}}^{\it ext}_g:={\mathcal{A}}\otimes g.\label{asg}
\end{equation}
In this particular case the $H$--Galois extension of $\mathcal{A}$ which is also the
$H$--Hopf module means that $\mathcal{A}^{\it ext}$ is $G$--graded algebra built from
a few copies (colors) of $\mathcal{A}$.

\section{Algebras of reality}

Let $H = {\bf k} G$, where $G$ is an Abelian group. We use here the concept of
$G$--graded $b$--commutative algebras and the so-called standard gradation \cite{WM4}.
This means that the grading group is $G\equiv Z^N := Z\oplus...\oplus Z$
($N$-sumands). In this case we have
\begin{equation}
b (\xi^i , \xi^j)=b^{ij} = (-1)^{\Sigma_{ij}} q^{\Omega_{ij}}, \label{comd}
\end{equation}
where $\{\xi^i : i = 1, \ldots , N\}$ is a set of generators of $G$, $\Sigma :=
(\Sigma_{ij})$ and $\Omega := (\Omega_{ij})$ are integer--valued matrices such that
$\Sigma_{ij} = \Sigma_{ji}$ and $\Omega_{ij} = - \Omega_{ji}$, $q \in {\bf k}
\setminus \{0\}$ is a parameter \cite{zoz}. We use here the following notation for the
generators of the grading group
\begin{equation}
\xi^i := (0, \ldots, 1, \ldots, 0),
\end{equation}
where $1$ is on the $i$--th place. If $q = exp(\frac{2 \pi i}{n})$, $n \underline{>}
3$, then the grading group $G \equiv Z^N$ can be reduced to $G = Z_n \oplus...\oplus
Z_n$. If $q = \pm 1$, then the grading group $G$ can be reduced to the group $Z_2
\oplus...\oplus Z_2$.

According to (\ref{asg}) we have $\mathcal{A}^{\it
ext}_{\xi^i}\simeq\mathcal{A}\otimes\xi_i$. We use the notation
\begin{equation}
x^a_{\xi^i}:=\theta^a\otimes\xi_i
\end{equation}
for generators of $\mathcal{A}^{\it ext}$. From (\ref{cod}), and (\ref{hco}) we obtain
\begin{equation}
\delta (x^a_{\xi^i}):= x^a_{\xi^i}\otimes\xi^i,\quad
x^a_{\xi^i}\triangleleft\xi^j:=x^a_{\xi^i\xi^j}\otimes \xi^i\xi^j.
\end{equation}
and we have the following quantum commutativity
\begin{equation}
x^a_{\xi^i}\;x^b_{\xi^j}:=b^{ij}\;x^b_{\xi^j}\; x^a_{\xi^i}.
\end{equation}
Assume that the algebra $\mathcal{A}$ is generated by $\theta^1,\ldots ,\theta^n$. We
define a mapping $r: \mathcal{A}^{\it ext}_{\xi^i}\rightarrow\mathcal{A}^{\otimes N}$
\begin{equation}
r(x^a_{\xi^i}):=(1\otimes\cdots\otimes\theta^a\otimes\cdots\otimes 1),
\end{equation}
where $\theta^a$ is on the $i$-th place. One can see that $r$ is a $H$-Hopf bimodule
morphism. If the mapping $r$ can be lifted to a well-defined algebra isomorphism,
then we say that ${\mathcal{A}}^{\it ext}$ is an algebra of reality describing
charges and quanta as a part of classical system.

Let us consider some simple examples. Let $\mathcal{A}$ be an algebra generated by one
Grassmann variable $\theta$, where $\theta^2 = 0$, and let the grading group $G$ be
the group $Z_2 \oplus...\oplus Z_2$. Then the algebra ${\mathcal{A}}^{\it ext}$ is
generated by $x^i$ and relations
\begin{equation}
x^i x^j = b^{ij} x^j x^i,\quad (x^i)^2 = 0 \quad \mbox{for}\quad b^{ii} = -1,
\end{equation}
where
\begin{equation}
\begin{array}{c}
b^{ij} = (-1)^{\Sigma^{ij}} (-1)^{\Omega^{ij}} \label{comf}
\end{array}
\end{equation}
The mapping $r$ is given by
\begin{equation}
r(x^i):=(1\otimes\cdots\otimes\theta\otimes\cdots\otimes 1).
\end{equation}
For electron in singular magnetic field we have \cite{top}
\begin{equation}
\begin{array}{c}
b^{ij} = (-1)^{\delta_{ij} + N} \label{com}
\end{array}
\end{equation}
where $N$ is the number of magnetic fluxes per particle.  If $N$ is even, then
$\epsilon^{ii} = -1$ for $i=1, 2, \ldots , N$, $\epsilon^{ij} = 1$ for $i\neq j$. For
$N=2$ we obtain
\begin{equation}
x^1 x^2 = x^2 x^1,
\end{equation}
and
\begin{equation}
r(x^1) := (\theta, 1),\quad r(x^2) := (1, \theta).
\end{equation}
The mapping $r$ is here an algebra isomorphism. In fact for $r(x^1 x^2)$ we obtain
\begin{equation}
r(x^1 x^2) = r(x^1)r(x^2) = (\theta, \theta). \label{haf}
\end{equation}
and similarly for $x^2 x^1$
\begin{equation}
r(x^2 x^1) = r(x^2) r(x^1) = (\theta, \theta). \label{paf}
\end{equation}
Note that the algebra $\mathcal{A}^{\it ext}$ is an example of the so--called $Z_2
\oplus Z_2$--graded commutative colour Lie superalgebra \cite{luri}. Observe that the
corresponding generators $x^1$ and $x^2$ of the algebra $\mathcal{A}^{\it ext}$
commute, their squares disappear and they describe two different quasiparticles. This
means that these two quasiparticles can become really existing particles. This also
means that single fermion $\theta$ can be transform under certain interactions into a
system of two different particles $x^1 x^2$. Such system is said to be a composite
fermion \cite{jai,sin}.

For $N$ odd we have $b^{ii} = 1$ for $i=1, 2, \ldots , N$, $b^{ij} = -1$ for $i\neq
j$. If $N=3$, then we obtain
\begin{equation}
\begin{array}{l}
r(x^1) =(\theta, 1, 1)\otimes\xi^1,\; r(x^2) = (1, \theta, 1)\otimes\xi^2,\; r(x^3) =
(1, 1, \theta)\otimes\xi^3. \label{sqa}
\end{array}
\end{equation}
Observe that in this case
\begin{equation}
x^i x^j = - x^j x^i ,\quad i \neq j .
\end{equation}
For the mapping $r$ we obtain $r(x^1 x^2) = (\theta, \theta, 1)$ and $r(x^2 x^1) =
(\theta, \theta, 1)$. But $r(x^2 x^1) = - r(x^1 x^2) = - (\theta, \theta, 1) = 0$.
Hence $x^1 x^2 = - x^2 x^1$ if and only if $x^1 x^2 = 0$, similar for $x^1 x^3$, $x^2
x^3$. This means that the multiplication in our extended algebra $\mathcal{A}^{\it
ext}$ is degenerate. This is just the generalized Pauli exclusion principle. In our
physical interpretation this means that the fermion $\theta$ can be transform into a
system $x^i$ called composite boson \cite{jai,sin}, where $i=1,2$ or $3$. This
situation correspond to fractional quantum Hall effect \cite{tsg,dst,sin}.


\begin{thebibliography}{99}
\bibitem{tsg} D. C. Tsui, H. L. Stormer, and A. C. Gossard,
{\em Phys. Rev. Lett.} {\bf 48} 1559 (1982)
\bibitem{eha} Z. F. Ezawa and H. Hotta, {\em Phys. Rev.} {\bf B 46},
7765 (1992)
\bibitem{hal} B. I. Halperin, P. A. Lee and N. Read,
{\em Phys. Rev.} {\bf B 46}, 7312 (1993)
\bibitem{zee} A. Zee, Quantum Hall fluids in Field Theory,
Topology and Condensed Matter Physics, ed. by H. D. Geyer, Lecture Notes in Physics,
Springer 1995.
\bibitem{jai} J. K. Jain, {\em Phys. Rev. Lett.} {\bf 63}, 199
(1989), {\em Phys. Rev.} {\bf B 40}, 8079 (1989); {\bf 41}, 7653 (1990)
\bibitem{dst} R. R. Du, H. L. Stormer, D. C. Tsui, A. S. Yeh,
L. N. Pfeiffer and K. W. West, {\em Phys. Rev. Lett.} {\bf 73}, 3274 (1994)
\bibitem{hal1} F. D. M. Haldane, {\it J. Phys.}
{\bf C14}, 2585 (1981).
\bibitem{bys} K. Byczuk and J. Spalek,
{\it Phys. Rev.} {\bf B51}, 7934 (1995).
\bibitem{qstat} W. Marcinek, Remarks on Quantum Statistics, in
Proceedings of the Conference "Particles, Fields and Gravitation", April 15 - 19,
(1998), Lodz, Poland, ed. by J. Rembielinski, World Scientific, Singapore 1998, and
math.QA./9806158.
\bibitem{com} W. Marcinek, On composite systems and quantum statistics,
in the Proceedings of the Vth International School on Theoretical Phsyics, Symmetry
and Structural Properties, Zajaczkowo k. Poznania, August 27 - September 2, 1998,
Poland, math.QA/9810060.
\bibitem{gqsd} W. Marcinek, On generalized quantum statistics,
in Proceedings of the XII-th Max Born Symposium, Wroc\l aw, September 23-26, 1998,
Poland.
\bibitem{gsd} W. Marcinek, On generalized statistics and one
dimensional systems, in Proceedings of the III International Seminar "Hidden
Symmetry", Rzeszow, October 20-22, 1998, Poland.
\bibitem{hor} L. P. Horwitz, {\em Found. Phys.} {\bf 22},
421 (1992).
\bibitem{haa} R. Haag, {\em Commun. Math. Phys.} {\bf 180},
733 (1996).
\bibitem{blj} ph. Blanchard, A. Jadczyk, {\em Ann. Phys.}
{\bf 4}, 583 (1995).
\bibitem{gsi} W. Marcinek, On generalized statistics and interactions,
in Proceedings of the XVI Workshop on Geometric Methods in Physics, July 1-7, 1998,
Bia{\l}owieza, Poland.
\bibitem{mon} S. Montgomery, Hopf algebras and their actions
on rings, Regional Conference series in Mathematics, {\bf No 82}, AMS 1993.
\bibitem{brmj} T. Brzezinski, S. Majid, Commun. Math. Phys. \textbf{191},
467-492, (1995).
\bibitem{sch} M. Scheunert, {\em J. Math. Phys.}
{\bf 20}, 712, (1979).
\bibitem{WM4} W. Marcinek, On unital braidings and quantization,
Rep. Math. Phys. 34, 325 (1994).
\bibitem{qsym} W. Marcinek, Particles and quantum symmetries,
in Proceedings of the XVI Workshop on Geometric Methods in Physics, July 1-7, 1997,
Bia{\l}owieza, Poland, mah.QA/9805122, to be pub. in Rep. Math. Phys.
\bibitem{zoz} Z. Oziewicz, Lie algebras for arbitrary grading group,
in Differential Geometry and Its Applications ed. by J. Janyska and D. Krupka, World
Scientific, Singapore 1990.
\bibitem{sma} S. Majid, Algebras and Hopf Algebras in Braided
Categories, in Advanced in Hopf Algebras, Plenum 1993.
\bibitem{bae} J. C. Baez, {\em Adv. Math.} {\bf 95}, 61 (1992).
\bibitem{mco} W. Marcinek, On algebraic model of composite fermions
and bosons, in  Proceedings of the IXth Max Born Symposium, Karpacz, September 25 -
September 28, 1996, Poland.
\bibitem{qweyl} W. Marcinek, On quantum Weyl algebras and generalized
quons, in Proceedings of the symposium:  Quantum Groups and Quantum Spaces, Warsaw,
November 20-29, 1995, Poland, ed. by R. Budzynski, W. Pusz and S. Zakrzewski, Banach
Center Publications, Warsaw 1997.
\bibitem{top} W. Marcinek, Topology and quantization, in Proceeding
of the IVth International School on Theoretical Physics, Symmetry and Structural
Properties, Zajaczkowo k. Poznania, August 29 - September 4 1996, Poland, hep-th/97
05 098, May 1997.
\bibitem{castat} W. Marcinek, Categories and quantum statistics,
in Proceedings of the symposium: Quantum Groups and their Applications in Physics,
Poznan October 17-20, 1995, Poland, {\it Rep. Math. Phys.} {\bf 38}, 149-179 (1996)
\bibitem{sin} W. Marcinek, {\em J. Math. Phys.}{\bf 39},
818--830 (1998).
\bibitem{wmq} W. Marcinek, {\em Rep. Math. Phys.} {\bf 41},
155 (1998).
\bibitem{cowe} M. Cohen and S. Westrich, {\em J. Alg.}
{\bf 168}, 1 (1994).
\bibitem{luri} J. Lukierski, V. Rittenberg,
{\it Phys. Rev.} {\bf D18}, 385, (1978).
\end{thebibliography}
\end{document}